\documentclass[preprint,11pt]{elsarticle}

\newtheorem{theorem}{Theorem}[section]
\newtheorem{lemma}{Lemma}[section]
\newtheorem{corollary}{Corollary}[section]

\newtheorem{remark}{Remark}[section]
\newcommand{\ignore}[1]{}{}

\usepackage{amssymb}
\usepackage{amsmath}

\usepackage{color}
\usepackage{soul}
\usepackage[dvipsnames]{xcolor}

\usepackage[colorlinks=true, linkcolor=blue, citecolor=blue]{hyperref}

\def\1{{{\mbox{${\rm{1\negthinspace\negthinspace I}}$}}}}

\newcommand\beq{\begin{equation}}
\newcommand\eeq{\end{equation}}

\usepackage{ifthen}
\usepackage{xkeyval}
\usepackage{todonotes}
\begin{document}

\begin{frontmatter}

\title{Self-normalized Cram\'{e}r  type moderate deviations for martingales and applications}
\author[cor1]{Xiequan Fan}
\author[cor2]{Qi-Man Shao}
\address[cor1]{School of Mathematics and Statistics, Northeastern University at Qinhuangdao, Qinhuangdao, 066004, Hebei, China.}
\address[cor2]{Department of Statistics and Data Science, SICM, National Center for Applied Mathematics Shenzhen,  Southern University of Science and Technology,
 Shenzhen  518000, China.}

\begin{abstract}
Cram\'er's moderate deviations give a quantitative estimate for the relative  error of the normal approximation
 and provide theoretical justifications for many estimator used in statistics. In this paper, we
 establish self-normalized Cram\'{e}r  type moderate deviations for martingales under some mile conditions. The result extends
 an earlier work of Fan, Grama, Liu and Shao [Bernoulli, 2019]. Moreover, applications of our result to Student's statistic,
 stationary martingale difference sequences and branching processes in a random environment are also discussed. In particular,
we establish Cram\'{e}r type  moderate deviations for Student's $t$-statistic for branching processes in a random environment.
\end{abstract}

\begin{keyword} Martingales; Self-normalized sequences; Cram\'{e}r's moderate deviations; Student's statistic; Branching process in a random environment
\vspace{0.3cm}
\MSC Primary 60G42; 60F10;  Secondary  60K37; 60J80
\end{keyword}

\end{frontmatter}




\section{Introduction}
\setcounter{equation}{0}

Let $(X_i)_{i\geq 1}$ be a sequence of independent and identically distributed (i.i.d.) random variables with means zero. Denote by  $\sigma= \sqrt{ \mathbf{E} X_1^2}$
 the standard   deviation of $X_1$. Assume that $\sigma > 0.$  Denote
$S_n=\sum_{i=1}^nX_i,  n\geq 1,$
the partial sum of $(X_i)_{i\geq 1}$.
The classical Cram\'{e}r moderate deviation \cite{Cr38} states that if $ \mathbf{E}e^{t_0 |X_1|} < \infty $ for some $t_0>0$, then
for all $0\leq x =o(\sqrt{n})$,
\begin{eqnarray} \label{fsdgsc1}
\bigg| \ln \frac{\mathbf{P}(  S_n/ \sqrt{n} \sigma \geq x)}{1-\Phi \left( x\right)} \bigg| \ = \   O\bigg( \frac{1+x^3}{ \sqrt{n}}  \bigg),\ \ \ \ n\rightarrow \infty.
\end{eqnarray}
In particular, the last equality implies that
\begin{eqnarray} \label{fsdgsc}
 \frac{\mathbf{P}(  S_n/ \sqrt{n} \sigma \geq x)}{1-\Phi \left( x\right)} \ = \ 1+ o(1)
\end{eqnarray}
holds uniformly for  $0\leq x =o(n^{1/6})$ as $n\rightarrow \infty$.  Roughly speaking, equality (\ref{fsdgsc}) shows that the relative error  for normal approximations on tail probability of standardized sum  tends to zero uniformly for $ 0\leq x =o(n^{1/6} )$. The results of types (\ref{fsdgsc1}) and (\ref{fsdgsc})
are usually called as Cram\'{e}r's  moderate deviations, and they have been well studied. We refer to  Chapter VIII of  Petrov \cite{Pe75} and references
therein for independent random variables  and  Saulis and Statulevi\v{c}ius \cite{SS78} for dependent random variables.
See also  Ra\v{c}kauskas \cite{Rackauskas90,Rackauskas95}, Grama and Haeusler \cite{GH00} and Fan \emph{et al.}\,\cite{FGL13} for martingales,
  Wu and Zhao \cite{WZ08} for stationary sequences with functional  dependence.

 One the other hand, the so-called self-normalized Cram\'{e}r type
moderate deviations put  a totaly new countenance upon the study for the relative errors of normal approximations.
Denote for $n\geq 1,$
$$  V_n^2= \sum_{i=1}^nX_i^2. $$
The self-normalized sum for independent random variables is defined as $S_n/V_n.$  Shao \cite{S97} gave  a self-normalized large deviation principle without any moment assumptions.
Shao \cite{S99} also has established the following self-normalized Cram\'{e}r type  moderate deviations: If   $\mathbf{E}|X_1|^{2+\rho}< \infty, \rho \in (0, 1] ,$ then it holds
\begin{equation}\label{Cramer01}
\frac{\mathbf{P}(S_n/ V_n \geq x)}{1-\Phi \left( x\right)}=1+o(1)
\end{equation}
uniformly for $0\leq x =o(n^{\rho/(4+2\rho)})$  as $n\rightarrow \infty$. Equality (\ref{Cramer01}) was further extended to independent but not necessarily identically distributed random variables by
 Jing, Shao and Wang  \citep{JSW03}. In particular, in the i.i.d. case, they
showed that it holds
\begin{equation}\label{Cramer02}
\bigg| \ln \frac{\mathbf{P}(S_n/ V_n \geq x)}{1-\Phi \left( x\right)} \bigg| = O\bigg( \frac{1+x^{2+\rho}}{ n^{\rho/2} }  \bigg)
\end{equation}
uniformly for $0\leq x =o(\sqrt{n})$. See also Liu, Shao and Wang \cite{LSW13} for self-normalized Cram\'{e}r type
moderate deviations for the maximum of sums $\max_{1\leq k \leq n}S_k/ V_n$ and  Shao and Zhou \cite{SZ16}  for general self-normalized processes.  We refer to de la Pe\~{n}a, Lai and Shao \cite{DLS09},  Shao and Wang \cite{SW13} and    Chang,   Shao and Zhou  \cite{CSZ16}
 for a systematic presentation on the self-normalized limit theory and its statistical applications.
Comparing to the classical Cram\'{e}r moderate deviations (\ref{fsdgsc1}) and (\ref{fsdgsc}),
self-normalized Cram\'{e}r moderate deviations of types  (\ref{Cramer01}) and (\ref{Cramer02}) have two significant  advantages in statistical inference of means, since
in practice one usually does not know the variance $\sigma^2.$
Moreover, self-normalized Cram\'{e}r type moderate deviations only require  finite moments
of the random variables instead of existence of exponential moments.

Despite the fact that the theory for self-normalized Cram\'{e}r type moderate deviations  for independent random variables
has been studied in depth, there are only a few results on expansions of types  (\ref{Cramer01}) and (\ref{Cramer02})  for   dependent
random variables. A few results in this direction, we refer to  Chen, Shao, Wu and Xu \cite{CSWX16},  Gao,  Shao  and Shi \cite{GSS22} for $\beta$-mixing and functional dependent sequences and Fan, Hu and Xu \cite{FHX23} for Euler-Maruyama scheme for SDE.
For some closely related topics, i.e.,  exponential inequalities for self-normalized martingales,
we refer to Bercu, Delyon and Rio \cite{BDR15} and de la Pe\~{n}a, Lai and Shao \cite{DLS09}.
Recently, Fan, Grama, Liu and Shao \cite{FGLS19} obtained the following self-normalized Cram\'{e}r type moderate deviations for  martingales.
 Let  $(B_n)_{n\geq 1}, (\delta_n)_{n\geq1}$   and $(\gamma_n)_{n\geq1}$  be three sequences of positive numbers, such that
$ B_n^{-1}, \delta_n, \gamma_n\rightarrow 0$ as $n\rightarrow \infty.$
Let $(X_i,\mathcal{F}_i)_{i\geq1}$ be a sequence of martingale differences satisfying
\begin{eqnarray}\label{fd78sd}
\Big|  \sum_{i=1}^n \mathbf{E} [ X _i^2  |  \mathcal{F}_{i-1} ] -B_n^2  \Big| \leq \delta_n^2 B_n^2\
\end{eqnarray} and
\begin{eqnarray}\label{fdsfsd}
 \mathbf{E}[|X_{i}| ^{2+\rho}  | \mathcal{F}_{i-1}]   \leq  (\gamma_{n }  B_n)^{\rho}\, \mathbf{E}[ X_{i}  ^{2}  | \mathcal{F}_{i-1}],
 \end{eqnarray}
 where $\rho \in (0, 1]$. Here and hereafter, the   inequalities between random variables are understood in the $\mathbf{P}$-almost sure sense.
 For $\gamma_n \rightarrow0,$ denote
 $$\widehat{\gamma}_n(x, \rho) =  \frac{ \gamma_n^{ \rho(2-\rho)/4 } }{ 1+ x  ^{  \rho(2+\rho)/4 }},\ \ \ \ \ x\geq0\textrm{ and }\rho \in (0, 1].$$
 By Corollary 2.3 in  \cite{FGLS19},  we have the following self-normalized Cram\'{e}r type moderate deviations:
for all $0\leq x =o(\gamma_n^{-1})$,
\begin{eqnarray}\label{gfdfsd}
\bigg| \ln \frac{\mathbf{P}(S_n/V_n \geq x)}{1-\Phi \left( x\right)} \bigg| \leq c_{\rho} \bigg( x^{2+\rho}  \gamma_n^\rho+ x^2 \delta_n^2 +(1+x)\big(  \widehat{\gamma}_n(x, \rho) +\delta_n  \big) \bigg)   .
\end{eqnarray}
In particular,   the last   inequality implies that
\begin{eqnarray}\label{ghklsd5}
 \frac{\mathbf{P}(S_n/V_n \geq x)}{1-\Phi \left( x\right)}   =   1+ o(1)
\end{eqnarray}
uniformly for $\displaystyle 0\leq x =o(\min\{ \gamma_n^{-\rho/(2+\rho)} , \delta_n^{-1  } \} ) $  as $n \rightarrow \infty.$
Although condition (\ref{fd78sd})  is very restrictive,  it still can be satisfied in a lot of pertinent examples, see Fan, Grama, Liu and Shao \cite{FGLS20}
for the examples satisfying condition (\ref{fd78sd}).
However, there also exist some  examples  that condition (\ref{fd78sd}) does not be satisfied. For instance, considering Student's statistic for
 branching process in a random environment, the related martingale  for the average of the Lotka-Nagaev estimators satisfies the following condition: For all $x>0,$ it holds
\begin{eqnarray}\label{dfds}
\mathbf{P} \bigg( \Big| \sum_{i=1}^n\mathbf{E}[ X_{i}^2  | \mathcal{F}_{i-1}] - B_n^2 \Big| \geq  x  B_n^2 \bigg ) \leq C \exp\bigg\{ -   x   \delta_n^{-2} \bigg\}.
\end{eqnarray}
This observation calls us to establish self-normalized Cram\'{e}r type moderate deviations under the conditions (\ref{fdsfsd}) and (\ref{dfds}).
In this paper, we prove that under conditions (\ref{fdsfsd}) and (\ref{dfds}), we have  for all $0\leq x  =o (\min\{ \gamma_n^{-1} , \delta_n^{-1} \})$,
\begin{eqnarray}\label{fdsds16}
\!\! \bigg|\ln \frac{\mathbf{P}(S_n/V_n \geq x)}{1-\Phi \left( x\right)} \bigg|  \leq    c_{\rho} \bigg(\!  x^{2+\rho}  \gamma_n^\rho +  x^2 \delta_n     +   (1+x ) \big(  \widehat{\gamma}_n(x, \rho)    +\delta_n |\ln \delta_n|\big)\!\!\bigg).
\end{eqnarray}
Compared to \eqref{gfdfsd},  the range of validity of the last inequality depends on $\delta_n$.
In particular, inequality (\ref{fdsds16})  implies that
\begin{eqnarray} \label{fdsds178}
 \frac{\mathbf{P}(S_n/V_n \geq x)}{1-\Phi \left( x\right)}   =   1+ o(1)
\end{eqnarray}
holds uniformly for $\displaystyle 0\leq x =o(\min\{ \gamma_n^{-\rho/(2+\rho)} , \delta_n^{-1/2 } \} ) $  as $n \rightarrow \infty.$
 Applications of (\ref{fdsds178}) to  \mbox{Student's} statistic,
stationary martingale difference sequences
and branching processes in a random environment  are also discussed.
In particular,  we establish
a self-normalized Cram\'{e}r type moderate deviation for branching processes $(Z_i)_{i\geq 0}$ in an i.i.d.\ random environment.
For two integers $n_0, n >0,$ denote
$$\hat{m}_{n_0,n}=\frac{1}{n}\sum_{k=n_0}^{n_0+n-1} \frac{Z_{ k+1}}{Z_{ k}}$$
the estimator of the offspring mean $m$. Let  $N_{n_0,n}$ be
 Student's $t$-statistic for the random variables $( \frac{Z_{ k+1}}{Z_{ k}} -m)_{ n_0 \leq k \leq n_0+n-1},$
 that is
\begin{eqnarray}\nonumber
N_{n_0,n} &=&\frac{ \sqrt{n}  (\hat{m}_{n_0,n} -m )}{\sqrt{   \frac{1}{n-1} \sum_{k=n_0}^{n_0+n-1} ( \frac{Z_{ k+1}}{Z_{ k}} - \hat{m}_{n_0,n} )^2 } \ } .
\end{eqnarray}
With some mild conditions, we prove that
\begin{eqnarray}\nonumber
\frac{\mathbf{P}(N_{n_0,n} \geq x)}{1-\Phi(x)} =1+o(1)
\end{eqnarray}
holds uniformly for $0\leq  x = o\big( n^{\rho/(4+2\rho)}   \big)$ as $n  \rightarrow \infty$.
Using this result, we establish the confidence intervals for
the offspring mean $m$.

The paper is organized as follows. Our main result is stated and discussed in Section \ref{sec2}.
In Section \ref{sec3}, we give some applications of our result
 to Student's $t$-statistic,  stationary martingale difference sequences and branching processes in a random environment.
The proofs of self-normalized Cram\'{e}r type moderate deviations
for stationary martingale difference sequences are given in
Section \ref{station}.
In Sections \ref{sec7} and \ref{sec8}, we give the proofs of self-normalized and standardized Cram\'{e}r type moderate deviations
for branching processes  in an i.i.d.\ random environment.

Throughout the paper,  denote by $c$ and $ c_\alpha$
a generic positive  constant and a generic positive constant depending only on $\alpha,$  respectively.
Their values may vary from line to line.
We also denote $\theta $ as any value  satisfying $\left| \theta  \right| \leq 1$.

\section{Main result}\label{sec2}
\setcounter{equation}{0}

Let  $(X_i,\mathcal{F}_i)_{i\geq 0 } $  be a  sequence of   martingale differences defined on a
 probability space $(\Omega ,\mathcal{F},\mathbf{P})$,  where $X_0=0 $ and
 $\{\emptyset, \Omega\}=\mathcal{F}_0\subseteq ...\subseteq \mathcal{F}_n\subseteq
\mathcal{F}$ are increasing $\sigma$-fields. Denote by
\begin{equation}
S_{0}=0,\ \ \ \ \ \ \ S_k=\sum_{i=1}^k X_i,\ \ \quad k\geq1.  \label{xk}
\end{equation}
Then $S=(S_k,\mathcal{F}_k)_{k\geq 0}$ is a martingale. 
Let $[S]$ and $\langle S \rangle$ be respectively  the square bracket and the conditional variance  of the
martingale $S,$   that is, for $k\geq1,$
\begin{eqnarray*}
[S]_0=0,\ \ \ \ \  [S]_k=\sum_{i=1}^k X_i^2, \ \ \ \ \ \langle S\rangle_0=0 \ \ \ \ \textrm{and} \ \ \ \langle S \rangle_k=\sum_{i=1}^k \mathbf{E} [ X_i^2  |  \mathcal{F}_{i-1} ].
\end{eqnarray*}
 Let  $(B_n)_{n\geq 1}, (\delta_n)_{n\geq1}$   and $(\gamma_n)_{n\geq1}$  be three sequences of positive numbers.
 Throughout the paper,  assume that  $\max\{B_n^{-1}, \delta_n,    \gamma_n \} \rightarrow 0 $ as $n\rightarrow \infty$.
We  introduce the following conditions:
\begin{description}
\item[\textbf{(A1)}]  There exist $B_n$ and $\delta_n\in (0, \frac14]$   such that for all $x>0,$
$$\mathbf{P} \Big( \big| \langle S \rangle_n - B_n^2 \big| \geq  x  B_n^2 \Big ) \leq c\, \exp\Big\{ -   x   \delta_n^{-2} \Big\};$$
\item[\textbf{(A2)}]   There exist a constant $\rho \in (0, 1]$ and   $ \gamma_n \in (0, \frac14]  $ such that
for all $1\leq i \leq n,$
\[
\mathbf{E}[|X_{i}| ^{2+\rho}  | \mathcal{F}_{i-1}]   \leq  (\gamma_{n }  B_n)^{\rho}\, \mathbf{E}[ X_{i}  ^{2}  | \mathcal{F}_{i-1}].
\]
\end{description}
Note that when $(X_{i})_{i\geq 1}$ is a sequence of i.i.d.\ random variables with finite  moments
of order $2+\rho$,   conditions  (A1) and (A2)  are satisfied with $\delta_n=O(1/\sqrt{n})$ and $\gamma_n=O(1/\sqrt{n})$.

For any sequence of positive numbers $(\gamma_n)_{n\geq 1}$,  denote
 \begin{equation}
 \label{alpphahat001}
 \widehat{\gamma}_n(x, \rho) =  \frac{ \gamma_n^{ \rho(2-\rho)/4 } }{ 1+ x  ^{  \rho(2+\rho)/4 }},\ \ \ \ \ x\geq0\textrm{ and }\rho \in (0, 1].
 \end{equation}
 Clearly, it holds $\widehat{\gamma}_n(x, \rho) \rightarrow 0$ when $x\rightarrow \infty$ or $ \gamma_n\rightarrow 0.$

 Our   main  result is the following Cram\'{e}r type moderate deviations for the self-normalized  martingale
$$W_n= S_n \big/  \sqrt{[ S]_n}.$$

\begin{theorem}\label{th0}
Assume that conditions (A1) and  (A2)  are satisfied. Then it holds
\begin{eqnarray}\label{dfgs5sdf}
\bigg|\ln \frac{\mathbf{P}(W_n \geq x)}{1-\Phi \left( x\right)} \bigg|  \leq    c_{\rho} \bigg( x^{2+\rho}  \gamma_n^\rho +  x^2 \delta_n     +   (1+x ) \big(  \widehat{\gamma}_n(x, \rho)    +\delta_n |\ln \delta_n|\big)\bigg)   .
\end{eqnarray}
uniformly for $0\leq x  =o (\min\{ \gamma_n^{-1} , \delta_n^{-1} \})$. In particular,  the last inequality implies that
\begin{eqnarray}\label{ghkl45}
 \frac{\mathbf{P}(W_n \geq x)}{1-\Phi \left( x\right)} \ = \ 1+ o(1)
\end{eqnarray}
holds uniformly for $\displaystyle 0\leq x =o(\min\{ \gamma_n^{-\rho/(2+\rho)} ,\, \delta_n^{-1/2} \} ).$
Moreover, the inequality (\ref{dfgs5sdf}) remains valid   when $\displaystyle  \frac{\mathbf{P}(W_n \geq x)}{1-\Phi \left( x\right)}$ is replaced  by $\displaystyle \frac{\mathbf{P}(W_n \leq -x)}{ \Phi \left( -x\right)}$.
\end{theorem}

The detailed proof of  Theorem \ref{th0} is given in the supplement \cite{FS23b}.

Equality (\ref{ghkl45}) states that the relative error for normal approximation on tail probability of $W_n$ tends to zero uniformly for $ 0\leq x =o(\min\{ \gamma_n^{-\rho/(2+\rho)} , \, \delta_n^{-1/2} \} )$ as $n \rightarrow \infty$.

\begin{remark}\label{fdsdds}
Consider the following condition:
\begin{description}
\item[(\textbf{A1$'$})]  There exists   $ \delta_n \in (0, \frac14]  $ such that for all  $x>0,$
$$\mathbf{P}  \Big( \big| \langle S \rangle_n - B_n^2 \big| \geq  x  B_n^2 \Big )\ \leq \ c\, \exp\Big\{ -   x^2  \, \delta_n^{-2} \Big\}.$$
\end{description}
Clearly, when $0<x \leq 1,$ condition (A1) implies condition (A1$'$). In this sense, condition (A1) is more restrictive than condition (A1$'$).
When condition (A1) is
replaced by (A1$'$), following the proof of Theorem \ref{th0} but with $\delta_n(\lambda)  =c_0 \big( \lambda  \delta_n  + \delta_n \sqrt{|\ln \delta_n|}  \,\big)$,
we have for all $0\leq x  =o (\min\{ \gamma_n^{-1} , \, \delta_n^{-1} \})$,
\begin{eqnarray*}
\bigg|\ln \frac{\mathbf{P}(W_n \geq x)}{1-\Phi \left( x\right)} \bigg|
 \leq  c_{\rho} \bigg( x^{2+\rho}  \gamma_n^\rho +  x^3 \delta_n     +   (1+x ) \big(  \widehat{\gamma}_n(x, \rho)   +\delta_n |\ln \delta_n|\big)\bigg).
\end{eqnarray*}
The last inequality implies that (\ref{ghkl45})
holds uniformly for $\displaystyle 0\leq x =o(\min\{ \gamma_n^{-\rho/(2+\rho)}, \, $ $ \delta_n^{-1/3} \} ).$
\end{remark}

 The following corollary is a direct consequence of  Theorem \ref{th0}.  Such type result  for independent random variables can be found in Jing \textit{et al.}\ \cite{JSW03}.

\begin{corollary}\label{corollary01}
Assume that the conditions of Theorem \ref{th0} are satisfied. Then
 \begin{eqnarray*}
 \frac{ \ln \mathbf{P}(W_n \geq x)}{\ln (1-\Phi \left( x\right))} \ = \ 1+ o(1)  \ \ \ \ \ \ \ \textrm{and} \ \ \ \ \ \ \ \frac{ \ln \mathbf{P}(W_n \leq -x)}{\ln   \Phi \left( -x\right) } \ = \ 1+ o(1)
\end{eqnarray*}
  hold   uniformly for $0\leq x  =o (\min\{ \gamma_n^{-1} , \delta_n^{-1} \})$.
\end{corollary}

By Theorem \ref{th0}, with an argument similar to the proof of Corollary 2.2 in \cite{FGLS19}, we can obtain  the following moderate  deviation principle (MDP) for the self-normalized martingale  $W_n$ with the  conditions (A1) and  (A2).
\begin{corollary}\label{corollary02}
Assume that the conditions of Theorem \ref{th0} are satisfied.
Let $a_n$ be any sequence of real numbers satisfying $a_n \rightarrow \infty$ and $a_n \min\{ \gamma_n,  \delta_n\} \rightarrow 0$
as $n\rightarrow \infty$.  Then  for each Borel set $B$,
\begin{eqnarray}
- \inf_{x \in B^o}\frac{x^2}{2} \ &\leq &\ \liminf_{n\rightarrow \infty}\frac{1}{a_n^2}\ln \mathbf{P}\bigg(\frac{  W_n}{a_n }  \in B \bigg) \nonumber \\
 &\leq&\ \limsup_{n\rightarrow \infty}\frac{1}{a_n^2}\ln \mathbf{P}\bigg(\frac{ W_n}{a_n }    \in B \bigg)   \, \leq   \, - \inf_{x \in \overline{B}}\frac{x^2}{2} \, ,   \label{MDP}
\end{eqnarray}
where $B^o$ and $\overline{B}$ denote the interior and the closure of $B$, respectively.
\end{corollary}

By Theorem \ref{th0},   we can   obtain the following uniform self-normalized  Berry-Esseen bound:
\begin{eqnarray}
\sup_{x \in \mathbf{R}}\Big|\mathbf{P}(W_n \leq x) -    \Phi(x) \Big| \leq  c_\rho    \Big(  \gamma_n^{ \rho(2-\rho)/4 }    +\delta_n |\ln \delta_n|\Big),
\end{eqnarray}
where $c_\rho$ depends only on  $\rho$. In particular, the last inequality implies the self-normalized central limit theory
 for $W_n$, that is, if $\gamma_n, \delta_n \rightarrow 0$ as $n\rightarrow \infty$, then $W_n$
converges to the standard normal random variable in distribution.

\section{Applications}\label{sec3}
\setcounter{equation}{0}

\subsection{Student's $t$-statistic}
The study of self-normalized  partial sums $S_n/V_n$
originates from Student's $t$-statistic.
 Student's $t$-statistic $T_n$ is defined as follows
\[
T_n = \sqrt{n} \, \overline{X}_n / \widehat{\sigma}_n,
\]
where $$\overline{X}_n = \frac{S_n}{n}  \ \ \ \ \ \ \ \textrm{and}\ \ \ \ \  \ \ \widehat{\sigma}_n^2 =\frac{1}{n-1} \sum_{i=1}^n  (X_i - \overline{X}_n )^2  .$$
The asymptotic distribution of $T_n$ is critical for the estimation of $p$-value for $t$-test.

When $(X_i)_{i\geq 1}$  are independent random variables,
Gin\'{e},  G\"otze and  Mason  \cite{GGM97} gave a necessary and sufficient condition for the asymptotic normality of $T_n$.
Slavova \citep{S85} and Bentkus and G\"{o}tze \cite{BG96} obtained Berry-Esseen's bounds for Student's $t$-statistic (see also \mbox{Novak} \cite{N11} for Berry-Esseen type inequalities with explicit constants).
Shao \cite{S99}  and Jing, Shao and Wang  \citep{JSW03} have established Cram\'{e}r type moderate deviations for $t$-statistic.

It is known that for all $x\geq0,$
\begin{equation} \label{chunf}
\mathbf{P}\Big( T_n  \geq x \Big) = \mathbf{P}\bigg(  W_n  \geq x \Big(\frac{n}{n+x^2-1} \Big)^{1/2}  \bigg ),
\end{equation}
see Chung \citep{C46}.
Thus Theorem \ref{th0} can also be applied to $t$-statistic for martingales.
Such type applications can be found in Fan \textit{et al.} \cite{FGLS19}.

\subsection{Stationary martingale difference sequences}

Assume that 
$X_i=X_0 \circ T^i $ and
$T: \Omega \mapsto \Omega$ is a bijective bimeasurable transformation preserving the probability $\mathbf{P}$ on $(\Omega, \mathcal{F})$.
For a subfield $\mathcal{F}_0$ satisfying $\mathcal{F}_0 \subseteq T^{-1}(\mathcal{F}_0)$, let $\mathcal{F}_i= T^{-i}(\mathcal{F}_0).$
Assume  that $X_0$ is   $\mathcal{F}_0$-measurable and  that $(X_i, \mathcal{F}_i)_{i \geq 0}$ is a stationary martingale difference sequence.
Denote by  $\lfloor a \rfloor$ the integer part of the real $a.$ Let
  $$m \in  [1, n] \ \ \ \ \ \ \textrm{and} \ \ \ \ \  \   k=\lfloor n/m   \rfloor,$$
where $m$ is an integer which may depend on $n.$ For instance, one may take $m=\lfloor n^c   \rfloor$ for a constant $c \in (0, 1).$
Define
$$H_j= \{i: m(j-1)+1  \leq i \leq mj  \}, \ \ \ \ 1 \leq j \leq k.$$
Denote by $S_{j}^\circ= \sum_{i \in H_j} X_i$ the block sum. The block self-normalized sum for $(X_i)_{1\leq i  \leq n}$ is defined by
$$ W^\circ_n= \frac{\sum_{j=1}^k S_{j}^\circ  }{ V_k^\circ   }, \ \    \ \ \ \ \ \ \textrm{where } \ \   (V_k^\circ)^2 =  \sum_{j=1}^k (S_{j}^\circ)^2   .$$
Clearly, $W^\circ_n$ is a self-normalized sum of $k$ random variables $(S_{j}^\circ)_{1\leq j \leq k}$.
We introduce the following conditions:
\begin{description}
\item[\textbf{(C1)}]  There exist a positive constant  $ \sigma$ and  $ \delta_m \in (0, \frac14] $ satisfying   $\delta_m^{2} \ln k \rightarrow 0$ as $ n\rightarrow \infty $  such that for all $x>0,$
\end{description}
\begin{equation} \label{gdfdfgdd}
\mathbf{P} \bigg( \Big  | \frac1{m\sigma^2}  \mathbf{E}[ S_m^2| \mathcal{F}_0]- 1\Big  | \geq x \bigg ) \leq c_1 \exp\bigg\{ -    x  \delta_m^{-2} \bigg\}.
\end{equation}
\begin{description}
\item[\textbf{(C2)}]   There exists a positive constant   $\rho \in (0, 1]$   such that
\[
\mathbf{E}[| S_m | ^{2+\rho}| \mathcal{F}_0] \, \leq c_2\, m^{  \rho/2} \, \mathbf{E}[  S_m^2 |\mathcal{F}_0].
\]
\end{description}

  In the following theorem, we give a self-normalized Cram\'{e}r type moderate deviation for stationary martingale difference sequences
under the conditions (C1) and (C2).
 For $n\geq 1$, we shall use below the notation  $\widehat{  \frac{1}{\sqrt{n}} }(x, \rho) , $
which means  sequences defined by (\ref{alpphahat001}) with $ \gamma_n$ replaced by $ \frac{1}{  \sqrt{n} }$.

 \begin{theorem}\label{thnds}
 Assume that the  conditions (C1) and (C2) are satisfied.
 Then it holds
\begin{eqnarray}\label{gdghcd}
  \bigg| \ln \frac{\mathbf{P}( W_n^\circ  \geq x   )}{1-\Phi \left(  x\right)}  \bigg|  \  &\leq&\  c_{\rho} \bigg( x^{2+\rho} \frac{1}{k^{ \rho/2}}   + x^2  \delta_m\sqrt{\ln k}    \\
&& \ \ \ \ \ \ \ \ \ +\ (1+x)\Big(  \delta_m\sqrt{\ln k} \Big| \ln (\delta_m\sqrt{\ln k})\Big| + \widehat{ \frac{1}{ \sqrt{k}} } (x, \rho)  \Big) \bigg)    \nonumber
\end{eqnarray}
uniformly for all $\displaystyle 0\leq x  = o\big(\min\big\{ \sqrt{k }  ,\  (\delta_m \sqrt{\ln k})^{-1 }   \big \} \big)$.     In particular, it  implies that
\begin{eqnarray}\label{thls}
 \frac{\mathbf{P}( W_n^\circ  \geq x    )}{1-\Phi \left(  x\right)}   = 1 +o(1) \
\end{eqnarray}
holds uniformly  for   $\displaystyle  0 \leq x = o \big( \min\big\{ k^{ \rho/(4+2\rho)} , \ (\delta_m \sqrt{\ln k})^{-1/2 } \big \} \big)$.
Moreover, the inequality (\ref{gdghcd}) remains valid   with $\displaystyle \frac{\mathbf{P}( W_n^\circ  \geq x  )}{1-\Phi \left(  x\right)}$ replacing by $\displaystyle \frac{\mathbf{P}(W_n^\circ \leq -x )}{ \Phi \left(-  x\right)}$.
\end{theorem}


Clearly, the inequality (\ref{gdfdfgdd}) is  implied by   the following condition
\begin{eqnarray} \label{cond02}
\delta_m = \Big \| \frac 1 {m\sigma^2}  \mathbf{E}[ S_m^2| \mathcal{F}_0]-1\Big \|_\infty .
\end{eqnarray}
 With the conditions (C2) and (\ref{cond02}), Fan \emph{et al.}\ \cite{FGLS20} gave some self-normalized Cram\'{e}r type moderate deviations for stationary processes.  Using Remark \ref{fdsdds}, we have the following improvement on Fan \emph{et al.}\ \cite{FGLS20} for  stationary martingale differences.
 \begin{theorem}\label{thndt}
 Assume that the  conditions (C2) and (\ref{cond02}) are satisfied for all $m\geq 1$.
 Then it holds
\begin{eqnarray}\label{gdsch2}
\bigg|\ln \frac{\mathbf{P}(W_n^\circ \geq x)}{1-\Phi \left( x\right)} \bigg|
 \leq  c_{\rho} \bigg( x^{2+\rho} \frac{1}{k^{ \rho/2}}+  x^3 \varepsilon_n    +   (1+x ) \Big(\varepsilon_n |\ln \varepsilon_n| +\widehat{ \frac{1}{ \sqrt{k}} } (x, \rho)   \Big)\bigg)
\end{eqnarray}
uniformly for  $0\leq x  =o (\min\{ \sqrt{k} , \, \varepsilon_n^{-1}  \})$, where $$\varepsilon_n=\frac{\sum_{j=1}^k j^{-1/2} \delta_{mj}}{\sqrt{k}} .$$
In particular, it  implies that
\begin{eqnarray}\label{dsfcsdg}
 \frac{\mathbf{P}( W_n^\circ  \geq x    )}{1-\Phi \left(  x\right)}   = 1 +o(1) \
\end{eqnarray}
holds uniformly  for   $\displaystyle  0 \leq x = o \big( \min\{  k^{ \rho/(4+2\rho)}, \, \varepsilon_n^{-1/3}  \}  \big)$.
Moreover,  the inequality (\ref{gdsch2}) remains  valid    with  $\displaystyle \frac{\mathbf{P}( W_n^\circ  \geq x  )}{1-\Phi \left(  x\right)}$ replacing by $\displaystyle \frac{\mathbf{P}(W_n^\circ \leq -x )}{ \Phi \left(-  x\right)}$.
\end{theorem}

Assume that $\delta_{m} =O( \frac{1}{\sqrt{m}})$ as $m\rightarrow \infty$. Then we have
$$\varepsilon_n =O\Bigg( \frac{\sum_{j=1}^k j^{-1/2} (mj)^{-1/2}  }{\sqrt{k}} \Bigg)= O \bigg(  \frac{\ln n  }{\sqrt{n}} \bigg) .  $$
Thus the range of validity of (\ref{dsfcsdg}) can be very large. In particular, taking $m=1$, we deduce that
\begin{eqnarray*}
 \frac{\mathbf{P}( W_n   \geq x    )}{1-\Phi \left(  x\right)}   = 1 +o(1) \
\end{eqnarray*}
 holds uniformly  for   $\displaystyle  0 \leq x = o \big( \min\{  n^{ \rho/(4+2\rho)}, n^{1/6}/(\ln n )^{1/3} \}  \big)$.

Notice that in the i.i.d.\,case, $W_n^o$ is a self-normalized sum  of $k$ i.i.d.\,random variables $(S_{j}^\circ)_{1\leq j \leq k} $
and  condition (C1) is satisfied with $\delta_m=1/n.$
According to the main result of Shao \cite{S99} (see also Jing \emph{et al.}\ \cite{JSW03}),  self-normalized Cram\'{e}r type moderate deviation    (\ref{thls}) holds  uniformly for $0\leq x =o(k^{ \rho/(4+2\rho)})$.
Thus the ranges of validity for  (\ref{thls})  and (\ref{dsfcsdg}) coincide  with the case
 of block self-normalized sums of i.i.d.\ random variables.

\subsection{Branching processes in a random environment}\label{subsec3.3}

 Limit theorems for branching process in a random environment (BPRE), introduced by  Smith  and Wilkinson \cite{SW69},  have  motivated many interesting works.
 See, for instance,   Afanasyev \emph{et al.}\ \cite{ABV14}, Vatutin and Zheng \cite{VZ12}  and  Bansaye  and  Vatutin \cite{BV17}
on the survival probability and conditional limit theorems for  subcritical BPRE.
 For supercritical  BPRE,   a number of works focus on moderate and large deviations; see, for instance,
B\"{o}inghoff  and Kersting \cite{BK10}, Bansaye and  B\"{o}inghoff \cite{BB11},
  Huang and Liu \cite{HL12},  Nakashima \cite{N13},   B\"{o}inghoff \cite{B14}  and Grama \emph{et al.}\ \cite{GLE17}.

A BPRE can be described as follows.   Let $\xi=(\xi_0, \xi_1, ...) $ be a sequence of i.i.d.\ random variables, where $\xi_n$ stands for the random environment at generation $n.$  For each $n \in \mathbf{N}=\{0, 1,...\}$,  the realization of $\xi_n$ corresponds a probability law $\{ p_i(\xi_n): i \in  \mathbf{N}\}.$
 A branching process $(Z_n)_{n\geq 0}$ in the random environment $\xi$ is given by the following formula:
 \begin{equation}
 Z_0=1,\ \ \ \ \ \ Z_{n+1}= \sum_{i=1}^{Z_n} X_{n,i}, \ \ \  \ \ \ \  n \geq 0,
 \end{equation}
where $X_{n,i}$ is the offspring number of   the $i$-th individual in the generation $n.$
Moreover, given the environment $\xi,$ the random variables $  (X_{n,i})_{i\geq 1} $ are independent of each other with a common  distribution law
as follows:
  \begin{equation} \label{ssfs}
 \mathbf{P}(X_{n,i} =k | \xi   )     = p_k(\xi_n),\ \ \ \ \ \ \ k\in \mathbf{N},
\end{equation}
and they are also independent of $(Z_i)_{1\leq i \leq n}.$
Denote by $\mathbf{P}_\xi$ the quenched law, i.e.\ the conditional probability when the environment
$\xi$ is given, and by $\tau$  the law of the environment $\xi.$ Then $\mathbf{P}(dx, d\xi )=\mathbf{P}_\xi(dx)\tau(d\xi)$
is the total law of the process $(Z_n)_{n\geq0}$, called annealed law. In the sequel, the expectations with respect to the quenched law and annealed  law are denoted
by $\mathbf{E}_\xi$ and $\mathbf{E}$, respectively.

  Denote by $m$  the average offspring number  of an individual, usually termed the offspring mean.
 Clearly, it holds  $$m=\mathbf{E}Z_1=\mathbf{E} X_{n,i} =\sum_{k=1}^\infty  k \mathbf{E} p_k(\xi_{n,i}),\ \ \ \ i \in \mathbf{N}.$$
Thus, $m$ is not only the mean of $Z_1$, but also the mean of each individual $X_{n,i}$.
Denote
\begin{eqnarray*}
  m_n:=m_n(\xi_n)=\mathbf{E}_\xi X_{n,i}=\sum_{i=1}^{\infty} i \, p_i(\xi_n) ,\ \ n\geq 0, \ \ \ \ \ \ \ \ \ \ \ \ \  \\
 \mu= \mathbf{E} \ln m_0,\ \ \ \ \nu^2= \mathbf{E}(  \ln m_0 -\mu)^2 \ \ \ \   \textrm{and} \ \ \ \  \tau^2=\mathbf{E} ( Z_1-m_0  )^2.
\end{eqnarray*}
 Recall that $\xi=(\xi_0, \xi_1, ...) $ is   a sequence of i.i.d.\ random variables.
Then $(m_n)_{n\geq 1}$ is also a sequence of nonnegative i.i.d.\ random variables. The distribution of $m_n$ depends only on $\xi_n$ and
$m_n$ is independent of $(Z_i)_{1\leq i \leq n}$. Clearly, it holds $m= \mathbf{E}m_n.$  The process $\{Z_n, n\geq 0 \}$  is respectively  called supercritical, critical or subcritical according to $\mu >0$, $\mu =0$ or $\mu <0$. Therefore, $\mu$ is  called as the critical parameter.

In this paper,   assume
\begin{eqnarray}\label{defv1}
 p_0(\xi_0)=0  \ \ \ a.s.,
\end{eqnarray}
which means each individual has at least one offspring. Thus the process $\{Z_n, n\geq 0 \}$ is   supercritical.
 Denote by $v$  and $\sigma$ the standard deviations of $Z_1$ and $m_0$ respectively,  that is
\begin{eqnarray*}
  \upsilon^2=\mathbf{E} (Z_1-m)^2,\ \ \ \ \ \ \ \ \  \sigma^2= \mathbf{E} (m_0-m)^2.
\end{eqnarray*}
To avoid triviality,   assume that  $0 < v, \sigma  <\infty,$ which implies that
  $\mu, \nu^2$ and $\tau$ all are finite.
The assumption $\sigma  >0$  means that the random environment is not degenerate.
Moreover,   assume that
\begin{eqnarray}\label{defv3}
\mathbf{E} \frac{Z_1}{m_0} \ln^+ Z_1  < \infty.
\end{eqnarray}
Denote
$$V_n= \frac{Z_n }{\Pi_{i=0}^{n-1}m_i},\ \ \ \ \ \ \ n \geq 1. $$
The  conditions (\ref{defv1}) and  (\ref{defv3}) together  implies
that the limit
 $$V=\lim_{n\rightarrow \infty} V_n   \ \ \ \textrm{exits a.s.},$$
$V_n \rightarrow V$ in $\mathbf{L}^{1} $ and   $\mathbf{P}(V >0)=\mathbf{P}(\lim_{n \rightarrow \infty} Z_n = \infty)=\lim_{n \rightarrow \infty} \mathbf{P}(Z_n >0) =1$
(see   Athreya and Karlin \cite{AK71b} and   Tanny \cite{T88}).
We also need the following assumption: There exist  two constants $p > 1$ and $\eta_0 \in (0, 1)$ such that
\begin{eqnarray}\label{defv4}
\mathbf{E} (\theta_0(p))^{\eta_0} < +\infty, \ \ \ \ \ \
\textrm{where} \ \ \  \theta_0(p)= \frac{Z_1^p}{m_0^p} .
\end{eqnarray}
The last assumption  implies that the harmonic moment is of a positive order, that is,
there exists a positive constant $\alpha$ such that
$\mathbf{E} V^{-\alpha}   < \infty,$ see Lemma \ref{ts021}.


A  critical task in statistical inference of BPRE is to estimate  the offspring mean   $m.$ To this end,
  the Lotka-Nagaev    \cite{L39,N67}   estimator $\frac{Z_{ k+1}}{Z_{ k}}$  plays an important role.
Recall that the conditions (\ref{defv1}) and (\ref{defv3}) together implies that  $Z_{n} \rightarrow \infty$ a.s.\ with respect to the annealed law $\mathbf{P}$.
By the law of large number, it is obvious that
  $\frac{Z_{ k+1}}{Z_{ k}}  \rightarrow m_k,  Z_k\rightarrow \infty,$ a.s.\ with respect to $\mathbf{P}_\xi$, and  that
    $\frac{1}{n} \sum_{k=1}^{n} m_k \rightarrow m   $ a.s.\ with respect to $\mathbf{P}$.
  Thus, we have $\frac{1}{n} \sum_{k=1}^{n} \frac{Z_{ k+1}}{Z_{ k}}   \rightarrow m $ a.s.\ with respect to $\mathbf{P}$,
  which means $\frac{1}{n} \sum_{k=1}^{n} \frac{Z_{ k+1}}{Z_{ k}}$ is a consistent estimator of $m.$
  Denote the estimator by $\hat{m}_{n_0,n}$, that is
  $$\hat{m}_{n_0,n}=\frac{1}{n}\sum_{k=n_0}^{n_0+n-1} \frac{Z_{ k+1}}{Z_{ k}}.$$
  It is natural to consider the following self-normalized process:
  For any $n_0\geq 0$ and $n\geq2$,  define
\begin{eqnarray}\label{fdsds}
N_{n_0,n}  = \frac{ \sqrt{n}  (\hat{m}_{n_0,n} -m )}{\sqrt{   \frac{1}{n-1} \sum_{k=n_0}^{n_0+n-1} ( \frac{Z_{ k+1}}{Z_{ k}} - \hat{m}_{n_0,n} )^2 } \ } .
\end{eqnarray}
Here $n_0$ may depend on $n.$
Clearly, $N_{n_0,n}$ can be regarded as Student's $t$-statistic for the random variables $( \frac{Z_{ k+1}}{Z_{ k}} -m)_{ n_0 \leq k \leq n_0+n-1}.$

\begin{remark} In (\ref{fdsds}), it is usual to  take $n_0=0$. However, in the real-world applications,
 it may happen that we know a historical data $(Z_{k})_{  n_0 \leq k \leq n_0+n}$ for some $n_0\geq0$ and $n\geq 2,$ 
 but the data $(Z_{k})_{ 0 \leq k \leq n_0-1}$ is missing.
  In such a case  $N_{0,n}$ is no longer applicable to estimate $m$, while $N_{n_0,n}$ can do it.  Motivated by this practical problem, it would be better to consider the more general case $n_0\geq0$ instead of taking $n_0=0.$
\end{remark}

The following theorem gives  self-normalized Cram\'{e}r  type moderate deviations for BPRE.

\begin{theorem}\label{th000}
Assume that $\mathbf{E}| Z_1  -m_0   |^{2+\rho } + \mathbf{E}   |  m_0 -m   |^{2+\rho } < \infty$ for some  $  \rho \in (0,  1]$.
Then it holds uniformly for  $0\leq x = o(   \sqrt{n/\ln n} ),$
\begin{equation}\label{ine001}
\bigg|\ln \frac{\mathbf{P}(N_{n_0,n} \geq x)}{1-\Phi \left( x\right)} \bigg|  =  O \bigg(    \frac{ x^{2+\rho} }{n^{\rho/2}}    +   (1+x)\widehat{ \frac{1}{ \sqrt{n}} } (x, \rho)         \bigg).
\end{equation}
In particular, it implies that
\begin{eqnarray}\label{ddsfafs}
\frac{\mathbf{P}(N_{n_0,n} \geq x)}{1-\Phi(x)}  =  1+o(1)
\end{eqnarray}
holds uniformly for $0\leq  x = o\big( n^{\rho/(4+2\rho)}   \big)$ as $n  \rightarrow \infty$.
Moreover,  equality (\ref{ine001}) remains  valid when $ \displaystyle \frac{\mathbf{P}( N_{n_0,n} \geq x)}{1-\Phi \left( x\right)}$ is replaced  by $\displaystyle \frac{\mathbf{P}(N_{n_0,n} \leq -x)}{ \Phi \left( -x\right)}$.
\end{theorem}


Closely related result to (\ref{ine001}) can be found in Grama \emph{et al.}\ \cite{GLE17}.
In particular, from Grama \emph{et al.}\ \cite{GLE17}, we have for all $0 \leq x = o(\sqrt{n} )$,
\begin{equation}\label{crasdmer}
\Bigg| \ln \frac{\mathbf{P}\big( \frac{ \ln  Z_{n}  - n \mu \ }{  \nu \sqrt{n}} \geq x  \big)}{1-\Phi(x)} \Bigg| = O\bigg( \frac{  1+x^3    }{  \sqrt{n}  }\bigg ).
\end{equation}
Recall that $m= \mathbf{E} m_0$ and $\mu= \mathbf{E} \ln m_0$. Moreover, by Jensen's inequality, the following relation holds
$$ e^\mu= \exp\{\mathbf{E} \ln m_0\} \leq  \mathbf{E} m_0 =m .$$
Comparing to (\ref{crasdmer}),  equality (\ref{ine001}) deals with the  estimation of $m$, instead of $\mu$.
Because equality (\ref{ine001})  does  not need to know the exact value of $\nu,$
moderate deviation results of types (\ref{ine001}) and
(\ref{ddsfafs}) play an important role in statistical inference of $m$ since in practice one usually does
not know the variance $\nu^2$.

By Theorem  \ref{th000}, following an argument similar to the proof of Corollary \ref{corollary02}, we obtain the following  MDP  result for $N_{n_0,n}$.
\begin{corollary}
Assume that the conditions  of Theorem  \ref{th000} are satisfied.
Let $(a_n)_{n\geq1}$ be any sequence of real numbers satisfying $a_n \rightarrow \infty$ and $a_n   \sqrt{\ln n }/ \sqrt{ n}      \rightarrow 0$
as $n\rightarrow \infty$.  Then   for each Borel set $B$,
\begin{eqnarray}
- \inf_{x \in B^o}\frac{x^2}{2} \, &\leq& \,  \liminf_{n\rightarrow \infty}\frac{1}{a_n^2}\ln \mathbf{P}\bigg(  \frac{N_{n_0,n} }{ a_n  }     \in B \bigg) \nonumber \\
  & \leq& \, \limsup_{n\rightarrow \infty}\frac{1}{a_n^2}\ln \mathbf{P}\bigg(\frac{ N_{n_0,n} }{ a_n }  \in B \bigg) \, \leq \,  - \inf_{x \in \overline{B}}\frac{x^2}{2}   ,   \label{SMDP}
\end{eqnarray}
where $B^o$ and $\overline{B}$ denote the interior and the closure of $B$, respectively.
\end{corollary}

 For any $n_0 \geq 0$ and $n\geq2$,  define
\begin{eqnarray}
S_{n_0,n}  = \frac{ \sqrt{n}}{ \sigma   \ }  (\hat{m}_{n_0,n} -m ) .
\end{eqnarray}
We have  Cram\'{e}r's moderate deviations for $S_{n_0,n}$.

\begin{theorem}\label{th00s1}
Assume that for all $l\geq 2,$
\begin{equation}\label{fdsdds01}
\mathbf{E}| Z_1  -m_0   |^{l } \ \leq \ \frac{1}{2} \  l! \  c^{l-2}\, \mathbf{E}| Z_1  -m_0   |^{2 }
\end{equation}
and
\begin{equation}\label{fdsdds02}
\mathbf{E}| m_0 -m   |^{l } \ \leq \ \frac{1}{2} \    l! \ c^{l-2}\, \mathbf{E}| m_0-m   |^{2 } .
\end{equation}
Then it holds uniformly for  $0\leq x = o(   \sqrt{n/\ln n} ),$
\begin{equation} \label{ine00s1}
\bigg|\ln \frac{\mathbf{P}(S_{n_0,n} \geq x)}{1-\Phi \left( x\right)} \bigg| =  O \bigg(\frac{x^3}{\sqrt{n}}+   (1+x)\frac{\ln n}{\sqrt{n}}      \bigg).
\end{equation}
In particular, it implies that
\begin{eqnarray} \nonumber
\frac{\mathbf{P}(S_{n_0,n} \geq x)}{1-\Phi(x)} =1+o(1)
\end{eqnarray}
holds uniformly for $0\leq  x = o\big( n^{1/6}   \big)$.
Moreover,  equality (\ref{ine00s1}) remains  valid when $ \displaystyle \frac{\mathbf{P}( S_{n_0,n} \geq x)}{1-\Phi \left( x\right)}$ is replaced  by $\displaystyle \frac{\mathbf{P}(S_{n_0,n} \leq -x)}{ \Phi \left( -x\right)}$.
\end{theorem}

The conditions (\ref{fdsdds01}) and (\ref{fdsdds02}) are called as Bernstein's conditions,
 and they are equivalent to
Cram\'{e}r's condition $ \mathbf{E}e^{t_0 |Z_1  -m_0 |} + \mathbf{E}e^{t_0 |m_0-m |}  < \infty $ for some $t_0> 0$,
 see Fan \textit{et al.} \cite{FGL13}.

By Theorem \ref{th00s1},   we can obtain the following Berry-Esseen bound:
\begin{eqnarray}
\sup_{x \in \mathbf{R}} \Big| \mathbf{P}(S_{n_0,n} \leq x) -\Phi \left( x\right) \Big|  =O\Big( \frac{\ln n}{\sqrt{n}} \Big)   .
\end{eqnarray}
Moreover, the last inequality implies the central limit theory for $S_{n_0,n}$, that is, $S_{n_0,n}$ converges to the standard normal random variable in distribution as $n\rightarrow \infty$.

The  following  MDP  result  is a simple consequence of Theorem \ref{th00s1}.
\begin{corollary}
Assume  that the conditions  of Theorem  \ref{th00s1} are satisfied. Then
(\ref{SMDP}) remains  valid when $N_{n_0,n}$ is replaced  by $S_{n_0,n}$.
\end{corollary}

\section{Proofs of Theorems \ref{thnds} and \ref{thndt}}\label{station}
\setcounter{equation}{0}

\textbf{Proof of Theorem \ref{thnds}.}
Recall that $(X_i, \mathcal{F}_i)_{i \geq 0}$ is a  stationary  sequence of martingale differences.
Then $(S_{j}^\circ, \mathcal{F}_{jm} )_{  j \geq 1}$ is also a   stationary  sequence of martingale differences.
By stationarity and condition (C2), it follows that
\begin{eqnarray}\label{gfdfd}
    \mathbf{E}\bigg[ \Big|\frac{S_j^\circ}{\sqrt{mk}\sigma }\Big|^{2+\rho}   \ \bigg| \mathcal{F}_{(j-1)m} \bigg ]   \leq  \Big( \frac{c_2^{1/\rho} }{\sigma  }\frac{1}{ \sqrt{k}} \ \Big)^\rho   \mathbf{E}\bigg[ \Big(\frac{S_j^\circ}{\sqrt{mk}\sigma }\Big)^{2}   \ \bigg| \mathcal{F}_{(j-1)m} \bigg ] .
  \end{eqnarray}
Notice that for all $x>0,$
\begin{eqnarray}
  \mathbf{P}\bigg( \Big | \frac1{mk\sigma^2} \sum_{j=1}^k  \mathbf{E}[ (S_j^\circ)^2 | \mathcal{F}_{(j-1)m}  ] - 1  \Big |  \geq x \bigg ) \
&=& \mathbf{P}\bigg( \Big | \frac1{m \sigma^2} \sum_{j=1}^k  \mathbf{E}[ (S_j^\circ)^2 | \mathcal{F}_{(j-1)m}  ] - k  \Big |  \geq k x \bigg ) \nonumber  \\
&\leq&   \mathbf{P}\bigg( \sum_{j=1}^k  \Big | \frac1{m \sigma^2} \mathbf{E}[ (S_j^\circ)^2 | \mathcal{F}_{(j-1)m}  ] - 1  \Big |  \geq k x \bigg ) \nonumber  \\
 &\leq& \sum_{j=1}^k\mathbf{P}\bigg( \Big | \frac1{m \sigma^2}   \mathbf{E}[ (S_j^\circ)^2 | \mathcal{F}_{(j-1)m}  ] - 1  \Big |  \geq  x  \bigg ) \nonumber \\
  &=& k \mathbf{P}\bigg(  \Big | \frac1{m \sigma^2}    \mathbf{E}[  S_m^2 | \mathcal{F}_{0}  ] - 1  \Big |  \geq  x \bigg )   \nonumber  \\
   &\leq& k c_1 \exp\bigg\{- x  \delta_m^{-2} \bigg\},  \label{gssds}
\end{eqnarray}
where the last line follows by condition (C1).
It is easy to see that  if $k \geq 3$, then it holds for all $x \geq 2\delta_m^2 \ln k  ,$
\begin{eqnarray}
  k  \exp\bigg\{- x  \delta_m^{-2} \bigg\} \ &\leq&\    \exp\bigg\{\ln k - \frac12 x \delta_m^{-2}- \frac12 x  \delta_m^{-2} \bigg\}\  \leq \ \exp\bigg\{ - \frac12 x  \delta_m^{-2} \bigg\} \nonumber  \\
  &\leq&\  \exp\bigg\{-\frac{x  \delta_m^{-2}}{2\ln k}  \bigg\}.\nonumber
\end{eqnarray}
Clearly, if $0\leq x < 2\delta_m^2 \ln k ,$ then it holds
$$\exp\bigg\{-\frac{x  \delta_m^{-2}}{2\ln k}  \bigg\} \leq e .$$
Thus, from (\ref{gssds}), we have for  all $x>0,$
\begin{eqnarray}
\mathbf{P}\bigg( \Big | \frac1{mk\sigma^2} \sum_{j=1}^k  \mathbf{E}[ (S_j^\circ)^2 | \mathcal{F}_{(j-1)m}  ] - 1  \Big |  \geq x \bigg )
    \leq  e \,c_1 \exp\bigg\{-\frac{x \delta_m^{-2}}{2\ln k}  \bigg\}.  \label{varsn}
\end{eqnarray}
 Set $X_j=  S_j^\circ, 1\leq j \leq k.$
 Then, by (\ref{gfdfd}) and (\ref{varsn}), the conditions (A1) and (A2) are satisfied with $$n=mk,\ \ \ \ B_n^2=mk\sigma^2,\ \ \ \ \displaystyle  \gamma_n  = \frac{c_2^{1/\rho} }{\sigma  } \frac{1}{\sqrt{k}\,} \ \ \ \ \
\textrm{and }  \ \ \  \delta_n =  \delta_m \sqrt{2\ln k } .$$  The desired result follows by Theorem \ref{th0}.
This completes the proof of Theorem \ref{thnds}.

The following  exponential inequality of  Peligrad \textit{et al.}\ \cite{PUW07} (cf.\ Proposition 2 therein) plays an important role in the proof  of  Theorem \ref{thndt}.
\begin{lemma}\label{lemma1}
  Let $(X_i)_{i \geq 1}$ be a stationary sequence of random variables adapted to the filtration $(\mathcal{F}_i)_{i \geq 1}$.
  Then, for all $x \geq 0,$
  \begin{equation}
\mathbf{P}\bigg( \max_{1\leq i \leq n} |S_i| \geq x \bigg)  \leq 4 \sqrt{e} \exp  \Bigg\{ - \frac{x^2 }{ 2n (  \|X_1\|_\infty + 80 \sum_{j=1}^n j^{-3/2}\|\mathbf{E}[S_j|\mathcal{F}_0]\|_\infty )^2 } \Bigg\}.
\end{equation}
\end{lemma}

\textbf{Proof of Theorem \ref{thndt}.}
Recall that $(S_{j}^\circ, \mathcal{F}_{jm} )_{  j \geq 1}$ is   a   stationary  sequence of martingale differences.
Notice that for all $x>0,$
\begin{eqnarray}
  \mathbf{P}\bigg( \Big | \frac1{mk\sigma^2} \sum_{j=1}^k  \mathbf{E}[ (S_j^\circ)^2 | \mathcal{F}_{(j-1)m}  ] - 1  \Big |  \geq x \bigg ) \
&=& \mathbf{P}\bigg( \Big | \frac1{m \sigma^2} \sum_{j=1}^k  \mathbf{E}[ (S_j^\circ)^2 | \mathcal{F}_{(j-1)m}  ] - k  \Big |  \geq k x \bigg ) \nonumber .
\end{eqnarray}
By Lemma \ref{lemma1} and condition (\ref{cond02}), we have  for all $x>0,$
\begin{eqnarray}
&& \mathbf{P}\bigg( \Big | \frac1{mk\sigma^2} \sum_{j=1}^k  \mathbf{E}[ (S_j^\circ)^2 | \mathcal{F}_{(j-1)m}  ] - 1  \Big |  \geq x \bigg )  \ \,\nonumber \\
 &&\   \leq  4 \sqrt{e} \exp  \Bigg\{ - \frac{k^2 x^2 }{ 2k \big(  \delta_m  + 80 \sum_{j=1}^k j^{-3/2} \| \mathbf{E}[ \sum_{l=1}^j ( \frac1{m \sigma^2} \mathbf{E}[ (S_l^\circ)^2 | \mathcal{F}_{(l-1)m}  ] - 1 ) | \mathcal{F}_{0} ]\|_\infty  )^2   } \Bigg\} \nonumber \\
 &&\   =  4 \sqrt{e} \exp  \Bigg\{ - \frac{k  x^2 }{ 2  \big(  \delta_m  + 80 \sum_{j=1}^k j^{-3/2} \|  \sum_{l=1}^j ( \frac1{m \sigma^2} \mathbf{E}[ (S_l^\circ)^2 | \mathcal{F}_{0}  ] - 1 )   \|_\infty  )^2   } \Bigg\} \nonumber \\
  &&\   = 4 \sqrt{e} \exp  \Bigg\{ - \frac{k  x^2 }{ 2  \big(  \delta_m  + 80 \sum_{j=1}^k j^{-1/2} \|    \frac1{m j \sigma^2} \mathbf{E}[  S_{mj}^2 | \mathcal{F}_{0}  ] - 1    \|_\infty  )^2   } \Bigg\} \nonumber \\
  &&\   = 4 \sqrt{e} \exp  \Bigg\{ - \frac{k  x^2 }{ 2  \big(  \delta_m  + 80 \sum_{j=1}^k j^{-1/2} \delta_{mj}  )^2   } \Bigg\} \nonumber \\
   &&\  \leq 4 \sqrt{e} \exp  \Bigg\{ - \frac{c_3 \, k  x^2 }{  \big(   \sum_{j=1}^k j^{-1/2} \delta_{mj}  )^2   } \Bigg\} .
\end{eqnarray}
 Set $X_j=  S_j^\circ, 1\leq j \leq k.$
 Then, by (\ref{gfdfd}) and the last inequality, the conditions (A1$'$) and (A2) are satisfied with $$n=mk,\ \ \ \ B_n^2=mk\sigma^2,\ \ \ \ \displaystyle  \gamma_n  = \frac{c_2^{1/\rho} }{\sigma  } \frac{1}{\sqrt{k}\,} \ \ \ \ \
\textrm{and }  \ \ \  \delta_n =  \frac{  \, \sum_{j=1}^k j^{-1/2} \delta_{mj}  }{\sqrt{c_3\, k\, } }    .$$
 The desired result follows by Remark \ref{fdsdds}.
This completes the proof of Theorem \ref{thndt}.

\section{Proof of Theorem \ref{th000}}\label{sec7}
\setcounter{equation}{0}

The following lemma can be found  in Grama \emph{et al.}\ \cite{GLP20}. Notice that
the condition $\mathbf{E}   |  m_0 -m   |^{2+\rho } < \infty$ implies  that $  \mathbf{E}  m_0  ^{2+\rho } < \infty$.
\begin{lemma}  \label{ts021}
There exists a positive constant $\alpha$ such that
$$\mathbf{E} V^{-\alpha}   < \infty.$$
\end{lemma}

Denote   $c_\rho$ a positive constant depending   only on   $\rho,  \upsilon^2 , \sigma^2, \tau^2,$  $\mu$, $ \nu, \alpha,$  $\mathbf{E} V^{-\alpha},$ $\mathbf{E}| Z_1  -m_0   |^{2+\rho }$ and $\mathbf{E}   |  m_0 -m   |^{2+\rho }$, but it does not depend on $n$ and $x$. 
Moreover, the exact values of $c_\rho$ may vary from line to line.
Using Lemma \ref{ts021}, we can easily obtain the following lemma.
\begin{lemma}  \label{thsn20}
   It holds for all $0<x \leq \mu \sqrt{n} /\nu,$
$$\mathbf{P}\bigg( \Big| \frac{\ln Z_{n}  - \mu  n }{\sqrt{ n}\, \nu }\Big| \geq x \bigg) \leq  c_{\rho}^{-1} \exp\Big\{    - c_{\rho}\, x^2 \Big\}.$$
\end{lemma}

\textbf{Proof.}
Notice that for all  $x> 0,$
\begin{eqnarray}
\mathbf{P}\bigg( \Big| \frac{\ln Z_{n}  - \mu  n }{\sqrt{ n}\, \nu }\Big| \geq x \bigg) &=& I_1 +I_2,
\end{eqnarray}
where
$$I_1=  \mathbf{P}\bigg(   \frac{\ln Z_{n}  - \mu  n }{\sqrt{ n}\, \nu }  \leq - x \bigg) \ \ \ \ \ \  \textrm{and} \ \ \ \ \ \ \ \ I_2 = \mathbf{P}\bigg(   \frac{\ln Z_{n}  - \mu  n }{\sqrt{ n}\, \nu }  \geq x \bigg).$$
Clearly, the following decomposition holds:
\begin{equation}\label{decopos}
 \ln Z_n = \sum_{i=1}^n X_i + \ln V_n,
\end{equation}
where $X_i=\ln m_{i-1} (i\geq 1)$ are i.i.d.\ random variables depending only on the environment $\xi.$
Denote $\eta_{i}=(X_i-\mu)/\sqrt{n}\, \nu.$
Clearly, it holds for all  $  x >0,$
 \begin{eqnarray}
   I_1
  &=&    \mathbf{P}\bigg(  \sum_{i=1}^n \eta_{i} + \frac{\ln V_{n}}{  \sqrt{n}\, \nu} \leq - x \bigg )   \nonumber \\
   &\leq&   \mathbf{P}\bigg( \sum_{i=1}^n \eta_{i}  \leq - \frac x 2 \bigg ) + \mathbf{P}\bigg(  \frac{\ln V_{n}}{  \sqrt{n}\, \nu} \leq -  \frac x 2  \bigg ) . \label{dasa01}
\end{eqnarray}
The condition $\mathbf{E}   |  m_0 -m   |^{2+\rho } < \infty$ implies  Cram\'{e}r's condition, that is $  \mathbf{E} e^{c_\rho\, |X_1| }  < \infty$ for
some $c_\rho > 0$.
By Corollary 1.2 of Liu and Watbled \cite{LW09}, we have
\begin{displaymath}
  \mathbf{P}\bigg( \frac{1}{n}\sum_{i=1}^n X_i  -\mu \leq - x \bigg )  \leq \left\{ \begin{array}{ll}
\exp\Big\{ - c_\rho \,  n \, x^2    \Big\},\ \ \ \  & \textrm{ $0< x\leq1$,}\\
  &  \\
\exp\Big\{ - c_\rho \, n \, x    \Big\},\ \ \ \  & \textrm{ $x> 1$. }
\end{array} \right.
\end{displaymath}
Thus, it holds for all $ 0  <x \leq \mu \sqrt{n} /\nu, $
 \begin{eqnarray}
 \mathbf{P}\bigg( \sum_{i=1}^n \eta_{i}  \leq - \frac x 2 \bigg )
 \, \leq \, \exp\bigg\{ - c_\rho \,  x^2    \bigg\}.
\label{dasa02}
\end{eqnarray}
By Markov's inequality,    it is easy to see that for  all $x>0,$
 \begin{eqnarray}
\mathbf{P}\bigg(  \frac{\ln V_{n}}{  \sqrt{n}\, \nu} \leq -  \frac x 2  \bigg )  &=&\mathbf{P}\bigg(  V_{n}^{-\alpha  } \geq \exp\Big\{\frac x 2 \sqrt{n} \nu \alpha \Big\} \bigg) \nonumber  \\
&\leq &  \exp\Big\{ - \frac x 2 \sqrt{n} \nu \alpha  \Big \}  \mathbf{E}  V_{ n}  ^{-\alpha } .\nonumber
\end{eqnarray}
where $\alpha$ is given by Lemma \ref{ts021}.
As $V_n \rightarrow V$ in $\mathbf{L}^{1} ,$   we have $V_n=\mathbf{E}[V| \mathcal{F}_n]$ a.s.\ By Jensen's inequality,
we deduce that
 \begin{eqnarray*}
 V_{ n} ^{-\alpha } =(\mathbf{E}[V  | \mathcal{F}_{ n}])^{-\alpha  } \leq \mathbf{E}[  V  ^{-\alpha } | \mathcal{F}_{ n}].
\end{eqnarray*}
Taking expectations on both sides of the last inequality, we deduce that
  \[
\mathbf{E} V_{n} ^{-\alpha} \leq \mathbf{E}  V ^{-\alpha  }  .
\]
Hence, by Lemma \ref{ts021}, we have for all  $ 0  <x \leq \mu \sqrt{n} /\nu, $
 \begin{eqnarray}
\mathbf{P}\bigg(  \frac{\ln V_{n}}{  \sqrt{n}\, \nu} \leq -  \frac x 2  \bigg )
 &\leq &   c_\rho \, \exp\bigg\{ - \frac x 2 \sqrt{n} \nu \alpha  \bigg \}  \leq  c_\rho \,  \exp\bigg\{ - c_\rho '  x^2    \bigg\} .   \label{dasa03}
\end{eqnarray}
Combining (\ref{dasa01})-(\ref{dasa03}) together, we get   for  all $ 0  <x \leq \mu \sqrt{n} /\nu, $
 \begin{eqnarray}
I_1   \leq    c_\rho^{-1}  \exp\bigg\{ - c_\rho    x^2    \bigg\} .\nonumber
\end{eqnarray}
Similarly, we can prove that the last bound holds also for $I_2.$ This completes the
proof of   lemma.

Now, we are in position to prove Theorem \ref{th000}.
Denote
 \begin{eqnarray}\label{gdfbg}
  \tilde{\xi}_{k+1}=  \frac{Z_{ k+1}}{Z_{ k}}  -m  ,
 \end{eqnarray}
$\mathfrak{F}_{n_0} =\{ \emptyset, \Omega \}  $ and $\mathfrak{F}_{k+1}=\sigma \{ \xi_{i-1}, Z_{i}: n_0\leq i\leq k+1  \}$ for  all $k\geq n_0$.
Notice that $X_{k,i}$ is independent of $Z_k.$
Then it is easy to verify that
 \begin{eqnarray}
\mathbf{E}[ \tilde{\xi}_{k+1}  |\mathfrak{F}_{k } ] &=&  Z_{ k} ^{-1 }  \mathbf{E}[   Z_{ k+1}  -mZ_{ k}  |\mathfrak{F}_{k } ] = Z_{ k} ^{-1 } \sum_{i=1}^{Z_k} \mathbf{E}[   X_{ k, i}  -m   |\mathfrak{F}_{k } ]  \nonumber \\
& =&   Z_{ k} ^{-1 } \sum_{i=1}^{Z_k} \mathbf{E}[   X_{ k, i}  -m   ] \nonumber \\
& =&0.
\end{eqnarray}
 Thus
 $(\tilde{\xi}_k, \mathfrak{F}_k)_{k=n_0+1,...,n_0+n}$ is a finite sequence of martingale  differences.
 Moreover, the following   equalities hold
\begin{eqnarray}
 \sum_{k=n_0}^{n_0+n-1} \mathbf{E}[ \tilde{\xi}_{k+1}^2  |\mathfrak{F}_{k } ] &=&
\sum_{k=n_0}^{n_0+n-1} \!\! Z_{k}^{-2} \mathbf{E}[  ( Z_{ k+1} -mZ_{ k} )^2  |\mathfrak{F}_{k } ]\nonumber \\
 &=& \sum_{k=n_0}^{n_0+n-1} \!\! Z_{k}^{-2} \mathbf{E}\Big[  \Big( \sum_{i=1}^{Z_k}    (X_{ k, i}  -m)  \Big)^2  \Big|\mathfrak{F}_{k } \Big] \nonumber \\
&=&\sum_{k=n_0}^{n_0+n-1}\!\!  Z_{k}^{-2} \mathbf{E}\Big[  \Big( \sum_{i=1}^{Z_k}    (X_{ k, i} -m_k+ m_k -m)  \Big)^2  \Big|\mathfrak{F}_{k } \Big]  \nonumber  \\
&=&\sum_{k=n_0}^{n_0+n-1}\!\!  Z_{k}^{-2} \Bigg(\mathbf{E}\Big[  \Big( \sum_{i=1}^{Z_k}    (X_{ k, i} -m_k)    \Big)^2 \Big |\mathfrak{F}_{k } \Big ] +\, \mathbf{E}\Big[  \Big( \sum_{i=1}^{Z_k}    (m_k -m ) \Big)^2  \Big|\mathfrak{F}_{k } \Big] \nonumber\\
  &&\ \ \ \ \ \ \ \ \ \ \ \ \ \ \  +\, 2\  \mathbf{E}\Big[  \Big( \sum_{i=1}^{Z_k}    (X_{ k, i} -m_k)    \Big)\Big( \sum_{i=1}^{Z_k}    (m_k -m ) \Big)  \Big|\mathfrak{F}_{k } \Big]\Bigg).\nonumber
\end{eqnarray}
Notice that conditionally
on  $\xi_n$, the random variables $(X_{n,i})_{i\geq1}$ are i.i.d. Thus,  we can deduce that
\begin{eqnarray*}
\mathbf{E}\Big[  \Big( \sum_{i=1}^{Z_k}    (X_{ k, i} -m_k)    \Big)^2  \Big |\mathfrak{F}_{k }  \Big]&=& \mathbf{E}\Big[ \mathbf{E}\Big[  \Big( \sum_{i=1}^{Z_k}    (X_{ k, i} -m_k)    \Big)^2  \Big|\xi_k , \mathfrak{F}_{k } \Big]   \Big| \mathfrak{F}_{k } \Big] \\
&=& \mathbf{E}\Big[    \sum_{i=1}^{Z_k}   \mathbf{E}[ (X_{ k, i} -m_k)^2      | \xi_k , \mathfrak{F}_{k } ]   \Big| \mathfrak{F}_{k } \Big] \\
&=& \sum_{i=1}^{Z_k}    \mathbf{E}[ (X_{ k, i} -m_k)^2      |   \mathfrak{F}_{k } ]
\ =\ \sum_{i=1}^{Z_k}   \mathbf{E}  (X_{ k, i} -m_k)^2\\
&=& Z_k  \tau^2.
\end{eqnarray*}
Similarly, we have
\begin{eqnarray*}
\mathbf{E}\Big[  \Big( \sum_{i=1}^{Z_k}    (X_{ k, i} -m_k)    \Big)\Big( \sum_{i=1}^{Z_k}    (m_k -m ) \Big)  \Big| \mathfrak{F}_{k } \Big]
 =  0
\end{eqnarray*}
and
\begin{eqnarray*}
\mathbf{E}\Big[  \Big( \sum_{i=1}^{Z_k}    (m_k -m ) \Big)^2  \Big|\mathfrak{F}_{k } \Big] =   Z_k^2 \mathbf{E}     (m_k -m ) ^2= Z_k^2 \mathbf{E}     (m_0 -m ) ^2=  Z_k^2\sigma^2 .
\end{eqnarray*}
Therefore, we get
\begin{eqnarray}\label{dsfsa}
  \mathbf{E}[ \tilde{\xi}_{k+1}^2  |\mathfrak{F}_{k } ] \geq  \mathbf{E}     (m_0 -m ) ^2=  \sigma^2
\end{eqnarray}
and
\begin{eqnarray}
 \sum_{k=n_0}^{n_0+n-1} \mathbf{E}[ \tilde{\xi}_{k+1}^2  |\mathfrak{F}_{k } ] = n \sigma^2     + \tau^2 \sum_{k=n_0}^{n_0+n-1} Z_{k}^{-1} .  \label{ineq9.2}
\end{eqnarray}
 Notice that $Z_{n_0+k}\geq Z_k.$
When $x\geq \frac{10 \ln n}{n \mu} ,$ we have
\begin{eqnarray}
& &\ \ \ \ \ \mathbf{P}\bigg( \frac1n \sum_{k=n_0}^{n_0+n-1} Z_{k}^{-1} \geq x \bigg)\leq \mathbf{P}\bigg( \frac1n \sum_{k=0}^{ n-1} Z_{k}^{-1} \geq x  \bigg) \nonumber  \\
&& \ \ \ \ \  \leq \ \mathbf{P}\bigg( \frac1n \sum_{k=0}^{ n-1} Z_{k}^{-1} \geq x,\ Z_{ [\frac34 nx]} \leq e^{ \mu [\frac12 nx]} \bigg)  \  + \ \mathbf{P}\bigg( \frac1n \sum_{k=0}^{ n-1} Z_{k}^{-1} \geq x ,   \ Z_{ [\frac34 nx]} > e^{\mu [\frac12 nx]} \bigg)  \nonumber  \\
&&\ \ \ \ \  \leq \ \mathbf{P}\bigg(   Z_{ [\frac34 nx]} \leq e^{ \mu [\frac12 nx]} \bigg)+  \mathbf{P}\bigg( \frac1n \sum_{k=0}^{ n-1} Z_{k}^{-1} \geq x ,   \ Z_{ [\frac34 nx]} > e^{ \mu [\frac12 nx]} \bigg) \nonumber  \\
& &\ \ \ \ \  \leq   \ I_1 +I_2, \nonumber
\end{eqnarray}
where
$$I_1=\mathbf{P}\Bigg(  \frac{\ln Z_{ [\frac34 nx]}  - \mu [\frac34 nx] }{ \sqrt{ [\frac34 nx]} \nu } \leq \frac{ \mu [\frac12 nx]- \mu [\frac34 nx] }{\sqrt{ [\frac34 nx]} \nu } \Bigg)  \ \ \
\textrm{and} \ \ \
  I_2=\mathbf{P}\Bigg(  \frac{ [\frac34 nx]}{n}  + \frac{n}{e^{ \mu [\frac12 nx]} } \geq x  \Bigg).$$
For the term $I_1,$ by Lemma  \ref{thsn20}, we have  for all $0< x \leq 1,$
\begin{eqnarray*}
I_1 &\leq&c_{\rho}^{-1}  \exp\Bigg\{- c_\rho \bigg(\frac{ \mu [\frac12 nx]- \mu [\frac34 nx] }{\sqrt{ [\frac34 nx]} \nu }  \bigg)^2  \Bigg\} \\
&\leq& c_{\rho}^{-1} \exp\Big\{-c'_{\rho} n\, x    \Big \}.
\end{eqnarray*}
For the  term $I_2$, it holds for all  $x\geq \frac{10 \ln n}{n \mu}$ and $n \geq 2+ \frac{2}{5} \mu e^{\mu }$,
$$\frac{ [\frac34 nx]}{n}  + \frac{n}{e^{ \mu [\frac12 nx]} } \leq \frac34 x + \frac{n}{e^{\frac12   10  \ln n - \mu } }=\frac34 x + \frac{e^{\mu } }{ n^4}   <x,$$
which implies that $I_2=0.$
Thus, we get  for all $0< x \leq 1$ and $n \geq 2+ \frac{2}{5} \mu e^{\mu }$,
\begin{eqnarray*}
\mathbf{P}\bigg( \frac1n \sum_{k=n_0}^{n_0+n-1} Z_{k}^{-1} \geq x \bigg)  &\leq &   c_{\rho}^{-1} \exp\Big\{-c'_{\rho } n\, x   \Big \}\mathbf{1}_{\{x\geq \frac{10 \ln n}{n \mu} \}} +  \mathbf{1}_{\{0< x< \frac{10 \ln n}{n \mu} \}}\\
&\leq&  c_{\rho}^{-1} e^{ 10c_{\rho}/\mu  } \exp\bigg\{-c_{\rho}\, \frac{ n}{\ln n}\, x   \bigg\}  .
\end{eqnarray*}
As $\frac1n \sum_{k=n_0}^{n_0+n-1} Z_{k}^{-1}\leq 1,$ the last line also holds for all $x> 1.$
Therefore, from (\ref{ineq9.2}),  we have for  all $x>0,$
\begin{eqnarray}
\mathbf{P}\Bigg(\bigg| \frac{1}{n \sigma^2 } \sum_{k=n_0}^{n_0+n-1} \mathbf{E}[ \tilde{\xi}_{k+1}^2  |\mathfrak{F}_{k } ]-1 \bigg| \geq x\Bigg)\leq  c_{\rho}^{-1} e^{ 10c_{\rho}/\mu  }   \exp\bigg\{-c_{\rho}\, \frac{ n}{\ln n}\, x    \bigg\} . \label{ine435}
\end{eqnarray}
 Notice that
\begin{eqnarray}
  \mathbf{E}[ |\tilde{\xi}_{k+1}|^{2+\rho} |\mathfrak{F}_{k } ] &=&
  Z_{k}^{-2-\rho } \mathbf{E}[  | Z_{ k+1} -mZ_{ k} |^{2+\rho }  |\mathfrak{F}_{k } ]  \nonumber \\
& = &  Z_{k}^{-2-\rho }  \mathbf{E}\bigg[  \Big| \sum_{i=1}^{Z_k}    (X_{ k, i}  -m)  \Big|^{2+\rho }  \bigg|\mathfrak{F}_{k } \bigg]. \label{ines2.8}
\end{eqnarray}
By Minkowski's inequality, we have
\begin{eqnarray}
&& \mathbf{E}\bigg[  \Big| \sum_{i=1}^{Z_k}    (X_{ k, i}  -m)  \Big|^{2+\rho }  \bigg|\mathfrak{F}_{k } \bigg]
\leq  \Bigg(  \sum_{i=1}^{Z_k}  \Big( \mathbf{E}[  |  X_{ k, i}  -m   |^{2+\rho }  |\mathfrak{F}_{k } ]   \Big)^{\frac{1}{2+\rho} }  \Bigg)^{2+\rho} \nonumber \\
&&=\Bigg(  \sum_{i=1}^{Z_k}  \Big( \mathbf{E}[  |  X_{ k, i}  -m_k + m_k -m  |^{2+\rho }  |\mathfrak{F}_{k } ]   \Big)^{\frac{1}{2+\rho} }  \Bigg)^{2+\rho} \nonumber\\
&&\leq Z_k^{2+\rho}  2^{1+\rho}    \bigg( \mathbf{E}[  |  X_{ k, i}  -m_k   |^{2+\rho }  |\mathfrak{F}_{k } ] + \mathbf{E}[  |   m_k -m  |^{2+\rho }  |\mathfrak{F}_{k } ]  \bigg)   \nonumber \\
&& =   2^{1+\rho}  Z_k^{2+\rho}  \Big(      \mathbf{E}| Z_1  -m _0  |^{2+\rho } +  \mathbf{E}   |  m_0 -m   |^{2+\rho } \Big) . \label{dfsds4}
\end{eqnarray}
Notice that $\mathbf{E}[ \tilde{\xi}_{k+1} ^{2} |\mathfrak{F}_{k } ]\geq \sigma^2$  (cf.\,inequality (\ref{dsfsa})).
 From (\ref{ines2.8}), by (\ref{dfsds4}) and the last inequality,   we deduce that
\begin{eqnarray}
  \mathbf{E}[ |\tilde{\xi}_{k+1}|^{2+\rho} |\mathfrak{F}_{k } ] &\leq& 2^{1+\rho}  \Big(        \mathbf{E}| Z_1  -m_0   |^{2+\rho }  + \mathbf{E}   |  m_0 -m   |^{2+\rho } \Big )   \nonumber \\
& \leq & 2^{1+\rho}  \frac{ \mathbf{E}| Z_1  -m_0   |^{2+\rho }+  \mathbf{E}   |  m_0 -m   |^{2+\rho } }{ \mathbf{E} (m_0  -m)^2 }     \mathbf{E}[ \tilde{\xi}_{k+1}^{2} |\mathfrak{F}_{k } ]   . \label{ineq9.4}
\end{eqnarray}
Let $\eta_k =\tilde{\xi}_{n_0+k} $ and $\mathcal{F}_{k}=\mathfrak{F}_{n_0+k}$.
Thus, by  (\ref{ineq9.4}) and (\ref{ine435}),  $(\eta_k,  \mathcal{F}_{k})_{k=1,...,n}$ is a finite martingale difference sequences and
 satisfies the conditions  (A1) and (A2) with
$$\delta_n^2  \asymp \frac{\ln n }{ n}\ \ \  \textrm{and} \ \ \  \gamma_n\asymp \frac{1}{\sqrt{n}} .$$
Clearly, it holds
$$  \widetilde{N}_{n_0,n}=  \frac{ \sum_{k=1}^{n}  \eta_{k } }{\sqrt{\sum_{k=1}^{n} \eta_{k }^2 } \ }. $$
Applying Theorem \ref{th0} to $(\eta_k, \mathcal{F}_{k})_{k=1,...,n}$,   we obtain the following   equality: For all $0\leq  x =  o(  \sqrt{n/\ln n }),$
\begin{equation}\label{ineq020}
\bigg|\ln\frac{\mathbf{P}(\widetilde{N}_{n_0,n} \geq x)}{1-\Phi(x)} \bigg| =O \bigg(    \frac{ x^{2+\rho} }{n^{\rho/2}}   +x^2\sqrt{\frac{\ln n}{n}}     +     \frac{ 1+x  }{n^{\rho(2-\rho)/8} ( 1+ x  ^{  \rho(2+\rho)/4 })}       \bigg).
\end{equation}
Notice that
\begin{displaymath}
x^2\sqrt{\frac{\ln n}{n}}   \leq \left\{ \begin{array}{ll}
 \displaystyle  \frac{ x^{2+\rho} }{n^{\rho/2}} ,\ \ \ \  & \textrm{ $0  \leq x   \leq  \sqrt{ \ln n }$,}\\
  &  \\
 \displaystyle      \frac{ x}{n^{\rho(2-\rho)/8} ( 1+ x  ^{  \rho(2+\rho)/4 })}  ,\ \ \ \  & \textrm{ $  x  >  \sqrt{ \ln n }$. }
\end{array} \right.
\end{displaymath}
Thus, from \eqref{ineq020}, we have for all $0\leq  x =  o(  \sqrt{n/\ln n }),$
\begin{equation}\label{ineq022}
\bigg|\ln\frac{\mathbf{P}(\widetilde{N}_{n_0,n} \geq x)}{1-\Phi(x)} \bigg| =O \bigg(    \frac{ x^{2+\rho} }{n^{\rho/2}}    +      \frac{1+x }{n^{\rho(2-\rho)/8} ( 1+ x  ^{  \rho(2+\rho)/4 })}       \bigg).
\end{equation}
 Denote $$x_n= x \Big(\frac{n}{n+x^2-1} \Big)^{1/2}.$$
By equality \eqref{chunf},
it is easy to see that for all $x\geq0,$
\begin{equation}\label{dfs}
\mathbf{P}\Big(N_{n_0,n}  \geq x \Big) = \mathbf{P}\Big(  \widetilde{N}_{n_0,n}  \geq x_n \Big ).
\end{equation}
Clearly, it holds
\begin{equation}  \label{fdfsdf01}
\bigg|\ln\frac{\mathbf{P}( N _{n_0,n} \geq x)}{1-\Phi(x)} \bigg| \leq \bigg|\ln\frac{\mathbf{P}\big( \widetilde{N} _{n_0,n} \geq x_n \big)}{1-\Phi(x_n)} \bigg| + \bigg|\ln\frac{1-\Phi(x_n)}{1-\Phi(x)} \bigg|.
\end{equation}
Notice that for all  $0 \leq x =o(\sqrt{n})$, we have $ \big| x_n  -x  \big|  = O(\frac{x^2+1}{n}  )   $ and
\begin{equation} \label{fdfsdf02}
 \bigg|\ln\frac{1-\Phi(x_n)}{1-\Phi(x)} \bigg|  =O\Big( \frac{1+  x^3 }{n}  \Big) .
 \end{equation}
Applying (\ref{ineq022}), (\ref{dfs}) and (\ref{fdfsdf02})  to  the inequality  (\ref{fdfsdf01}), we obtain the desired  equality  for  all $0\leq x =o(  \sqrt{n/\ln n} ) .$

Notice that as $(-\eta_k, \mathcal{F}_{k})_{k=1,...,n}$ is also a finite sequence of martingale differences,  the argument for $(\eta_k, \mathcal{F}_{k})_{k=1,...,n}$ also holds for $(-\eta_k, \mathcal{F}_{k})_{k=1,...,n}$. Thus,   equality  \eqref{ine001}    remains
valid when $\displaystyle \frac{\mathbf{P}( N_{n_0,n}\geq x)}{1-\Phi \left( x\right)}$ is replaced  by $\displaystyle  \frac{\mathbf{P}(  N_{n_0,n} \leq -x)}{ \Phi \left( -x\right)}$.

\section{Proof of Theorem \ref{th00s1}}\label{sec8}
\setcounter{equation}{0}

To prove Theorem \ref{th00s1}, we should make use of  the following lemma for martingale differences $(X_i,\mathcal{F}_i)_{i\geq 0 } $.
The proof of this lemma is similar to the proof of Theorem 2.2 in \cite{FS23} but with $
  \delta_n(\lambda)  =c_0 \big( \lambda^2 \delta_n^2 + \delta_n^2 | \ln \delta_n  |   \,\big)$.
\begin{description}
\item[(B2)]  There exists $ \epsilon_n \in (0, \frac14]  ,  \epsilon_n\rightarrow 0 ,$  such that for all $1\leq i\leq n,$
\[
\mathbf{E}[|X_{i}| ^{k} | \mathcal{F}_{i-1}]     \leq \frac12\, k! \, (B_n \epsilon_n)^{k-2} \,\mathbf{E}[ X_{i}  ^{2}  | \mathcal{F}_{i-1}],\ \ \ \ \ \  k\geq 2  .
\]
\end{description}
\begin{lemma}\label{thxz2s}
Assume that the conditions (A1) and (B2) are satisfied.  Then for all $0 \leq  x   =o( \min\{\epsilon_n^{-1}, $ $ \delta_n ^{-1 } \}),$
it holds
\begin{equation}\label{t0ide1}
 \bigg| \ln \frac{\mathbf{P}(S_n >xB_n )}{1-\Phi \left( x\right)} \bigg| \ \leq \ c_{ p}   \bigg( x^3  \epsilon_n  + x^2 \delta_n +   (1+  x )  \big(  \delta_n|\ln \delta_n| + \epsilon_n|\ln \epsilon_n| \big) \bigg) .
\end{equation}%
\end{lemma}

Now, we are in position to prove Theorem \ref{th00s1}.
Recall the definition of $ \tilde{\xi}_{k+1} $ (cf.\ equality \eqref{gdfbg}) and $\mathfrak{F}_{k } $.
We can deduce  that for all $l\geq 2,$
\begin{eqnarray}
\mathbf{E}[ |\tilde{\xi}_{k+1}|^{l} |\mathfrak{F}_{k } ]&=&
  Z_{k}^{-l } \mathbf{E}[  | Z_{ k+1} -mZ_{ k} |^{l }  |\mathfrak{F}_{k } ]   =    Z_{k}^{-l }  \mathbf{E}\bigg[  \Big| \sum_{i=1}^{Z_k}    (X_{ k, i}  -m)  \Big|^{l }  \bigg|\mathfrak{F}_{k } \bigg]  \nonumber \\
&  =&  Z_{k}^{-l } \mathbf{E}\bigg[  \Big| \sum_{i=1}^{Z_k}    (X_{ k, i}-m_k) + \sum_{i=1}^{Z_k}(m_k  -m)  \Big|^{l }  \bigg|\mathfrak{F}_{k } \bigg] \nonumber \\
&  \leq&  2^{l-1}Z_{k}^{-l } \Bigg( \mathbf{E}\bigg[  \Big| \sum_{i=1}^{Z_k}    (X_{ k, i}-m_k)\Big|^{l } \bigg|\mathfrak{F}_{k } \bigg]  +  \mathbf{E}\bigg[   \Big|\sum_{i=1}^{Z_k}(m_k  -m)  \Big|^{l }  \bigg|\mathfrak{F}_{k } \bigg] \Bigg)   \nonumber.
\end{eqnarray}
By Minkowski's inequality, we have  for all $l\geq 2,$
\begin{eqnarray}
\mathbf{E}[ |\tilde{\xi}_{k+1}|^{l} |\mathfrak{F}_{k } ]
&  \leq&  2^{l-1}Z_{k}^{-l } \Bigg(\bigg(  \sum_{i=1}^{Z_k} (\mathbf{E}[  | X_{ k, i}-m_k  |^{l } \big|\mathfrak{F}_{k }])^{1/l} \bigg)^l  +  \bigg(  \sum_{i=1}^{Z_k} (\mathbf{E}[  | m_k  -m  |^{l } \big|\mathfrak{F}_{k }])^{1/l} \bigg)^l \ \Bigg) \nonumber  \\
&  \leq&  2^{l-1}   \Big( \mathbf{E}   | X_{ k, i}-m_k  |^{l }     +   \mathbf{E}  | m_k  -m  |^{l }   \Big)
  \nonumber \\
&  =&  2^{l-1}   \Big(   \mathbf{E}   \big| Z_1-m_0 \big|^{l }    +   \mathbf{E}   \big| m_0  -m   \big|^{l }  \Big). \nonumber
\end{eqnarray}
Notice that $\mathbf{E}[ \tilde{\xi}_{k+1} ^{2} |\mathfrak{F}_{k } ]\geq \sigma^2$  (cf.\,inequality (\ref{dsfsa})).
Thus, by the conditions (\ref{fdsdds01}) and (\ref{fdsdds02}), we get for all $l\geq 2,$
\begin{eqnarray}
  \mathbf{E}[ |\tilde{\xi}_{k+1}|^{l} |\mathfrak{F}_{k } ] &\leq &   2^{l-1} \frac12 \  l! \  c^{l-2} \   \Big(   \mathbf{E}   \big( Z_1-m_0 \big)^{2 }    +   \mathbf{E}   \big( m_0  -m   \big)^{2 }  \Big) \nonumber \\
  &\leq& \frac12 \, l! \,    (2c)^{l-2} \, 2C_\sigma \, \mathbf{E}[ \tilde{\xi}_{k+1} ^{2} |\mathfrak{F}_{k } ]     \nonumber\\
  &\leq& \frac12 \, l! \,     C ^{l-2} \, \mathbf{E}[ \tilde{\xi}_{k+1} ^{2} |\mathfrak{F}_{k } ],  \nonumber
\end{eqnarray}
where
$$C_\sigma=\frac{1}{\sigma^2} \Big( \mathbf{E}   \big( Z_1-m_0 \big)^{2 }    +   \mathbf{E}   \big( m_0  -m   \big)^{2 } \Big)  \ \ \ \ \  \textrm{and} \ \ \ \ \ C=1+2c+2C_\sigma.$$
By (\ref{ine435}), it holds
\begin{eqnarray}
\mathbf{P}\Bigg(\bigg| \frac{1}{n \sigma^2 }\!\! \sum_{k=n_0}^{n_0+n-1} \mathbf{E}[ \tilde{\xi}_{k+1}^2  |\mathfrak{F}_{k } ]-1 \bigg| \geq x\Bigg)\leq  c_{\rho}^{-1} e^{ 10c_{\rho}/\mu  }   \exp\bigg\{-c_{\rho}\, \frac{ n}{\ln n}\, x    \bigg\} .  \nonumber
\end{eqnarray}
Let $\eta_k =\tilde{\xi}_{n_0+k}   $ and $\mathcal{F}_{k}=\mathfrak{F}_{n_0+k}$.
Then  $(\eta_k,  \mathcal{F}_{k})_{k=1,...,n}$ is a finite martingale difference sequences and
 satisfies   conditions  (A1) and (B2) with
$$\delta_n^2  \asymp \frac{\ln n }{ n}\ \ \  \textrm{and} \ \ \  \gamma_n\asymp \frac{1}{\sqrt{n}} .$$
Clearly, it holds
$ S_{n_0,n}=  \frac{1}{\sigma\sqrt{n}\ } \sum_{k=1}^{n}  \eta_{k }. $
Applying  Lemma \ref{thxz2s} to $(\eta_k, \mathcal{F}_{k})_{k=1,...,n}$,   we obtain the following  equality: For all $0\leq  x =  o(  \sqrt{n/\ln n }),$
\begin{equation}\label{ineq1022}
\bigg|\ln\frac{\mathbf{P}(S_{n_0,n} \geq x)}{1-\Phi(x)} \bigg| =O \bigg(  \frac{x^3  }{\sqrt{n}}  + x^2\sqrt{\frac{\ln n}{n}}    + (1+x)\frac{\ln n}{\sqrt{n}}      \bigg).
\end{equation}
Notice that
\begin{displaymath}
x^2\sqrt{\frac{\ln n}{n}}  \ \leq \ \left\{ \begin{array}{ll}
 \displaystyle x  \frac{\ln n}{n},\ \ \ \  & \textrm{ $0  \leq x   \leq  \sqrt{ \ln n }$,}\\
  &  \\
 \displaystyle \frac{x^3  }{\sqrt{n}},\ \ \ \  & \textrm{ $  x  >  \sqrt{ \ln n }$. }
\end{array} \right.
\end{displaymath}
Thus, by \eqref{ineq1022}, we have for all $0\leq  x =  o(  \sqrt{n/\ln n }),$
\begin{equation*}
\bigg|\ln\frac{\mathbf{P}(S_{n_0,n} \geq x)}{1-\Phi(x)} \bigg| =O \bigg(\frac{x^3 }{\sqrt{n}}      + (1+x)\frac{\ln n}{\sqrt{n}}      \bigg),
\end{equation*}
which gives the first desired equality.

The same argument applies to  $(-\eta_k, \mathcal{F}_{k})_{k=1,...,n}$, we find that equality  (\ref{ine00s1})    remains
valid when $\displaystyle \frac{\mathbf{P}( S_{n_0,n} \geq x)}{1-\Phi \left( x\right)}$ is replaced  by $\displaystyle  \frac{\mathbf{P}(  S_{n_0,n}  \leq -x)}{ \Phi \left( -x\right)}$.

\section*{Acknowledgements}
%
Fan would like to thank Quansheng Liu for his helpful discussion on the harmonic moments for branching processes in a random environment.
Fan X. was partially supported by the National Natural Science Foundation
of China (Grant Nos.\,11971063 and 12371155). Shao Q.M. was partially supported by the National Natural Science Foundation
of China (Grant No.\,12031005) and Shenzhen Outstanding Talents Training Fund.

\end{document}